%this file contains the paper by Dominique Foata and Doron
%Zeilberger, entitled ``Combinatorial Proofs of Capelli's and
%Turnbull's Identities from Classical Invariant Theory."
%The file is to be compiled with TeX and plain with no other
%format. There will be 10 printed pages to come.

\magnification=1200
\hsize=11.25cm
\hoffset1.2cm
\def\sgn{\mathop{\rm sgn}}
\def\qed{\quad\raise -2pt\hbox{\vrule
     \vbox to 10pt{\hrule width 4pt \vfill\hrule}\vrule}}

\centerline{\bf COMBINATORIAL PROOFS OF CAPELLI'S}
\vskip 5pt
\centerline{\bf AND TURNBULL'S IDENTITIES FROM}
\vskip 5pt
\centerline{\bf  CLASSICAL INVARIANT THEORY} 
\vskip 2.8 mm
\centerline{\sevenrm BY}
\vskip 2.8mm
\centerline{Dominique
FOATA\footnote{$^*$}{D\'epartement de math\'ematique,
Universit\'e Louis-Pasteur, 7, rue Ren\'e Descartes,
F-67084 Strasbourg Cedex, France (foata@math.u-strasbg.fr).} and Doron
ZEILBERGER\footnote{$^{**}$}{Supported in part by NSF
grant DM8800663; \hfil\break
department of Mathematics, Temple
University, Philadelphia, PA 19122, U.S.A. 
(zeilberg@euclid.math.temple.edu).}} 
\vskip 24 pt plus 6pt
{\bf 0. Introduction.}
Capelli's [C] identity plays a prominent role in Weyl's [W]
approach to Classical Invariant Theory. Capelli's identity
was recently considered by Howe~[H] and Howe and
Umeda~[H-U]. Howe~[H] gave an insightful
representation-theoretic proof of Capelli's identity, and
a similar approach was used in [H-U] to prove
Turnbull's~[T] symmetric analog, as well as a new
anti-symmetric analog, that was discovered
independently by Kostant and Sahi~[K-S]. The Capelli,
Turnbulll, and Howe-Umeda-Kostant-Sahi identities
immediately imply, and were inspired by, identities of
Cayley (see~[T1]), Garding~[G], and Shimura~[S],
respectively.

In this paper, we give short combinatorial proofs of
Capelli's and Turnbull's identities, and raise the hope
that someone else will use our approach to prove the
new Howe-Umeda-Kostant-Sahi identity.

\medskip
{\bf 1. The Capelli Identity.}
Throughout this paper
$x_{i,j}$ are mutually commuting indeterminates
(``positions"), as are $o_{i,j}$ (``momenta"), and they
interact with each other via the ``uncertainty principle"
$$
p_{ij}x_{ij}-x_{ij}p_{ij}=h,
$$
and otherwise $x_{i,j}$ commutes with all the $p_{k,l}$
if $(i,j)\not=(k,l)$. Of course, one can take
$p_{i,j}:=h(\partial/\partial x_{i,j})$. Set
$X=(x_{ij})$, $P=(p_{ij})$ $(1\le i,j\le n)$.

\proclaim Capelli's Identity. For each positive integer~$n$ and
for $1\le i,j\le n$ let 
$$
\displaylines{
\rlap{$(1.1)$}\hfill
A_{ij}=\sum_{k=1}^n
x_{ki}p_{kj} + h(n-i) \delta _{ij}.\hfill\cr
\noalign{\hbox{Then}}
\rlap{\sevenrm (CAP)}\hfill
\sum_{\sigma \in {\cal S}_n}
\sgn(\sigma )\,A_{\sigma 1,1}\ldots A_{\sigma n,n}=\det
X.\det P.\hfill\cr}
$$

{\it Remark $1$}.
The Capelli identity can be viewed as a ``quantum analog"
ot the classical Cauchy-Binet identity
$\det XP=\det X.\det P$, when the entries of $X$ and $P$
commute, and indeed reduces to it when $h=0$. The
matrix $A$ is $X^t\,P$, with ``quantum correction"
$h(n-i)\delta _{i,j}$.

\medbreak
{\it Remark $2$}. 
Note that since not all indeterminates
commute, it is necessary to define order in the
definition of the determinant of~$A$. It turns out that
the determinant is to be evaluated by ``column
expansion" rather than ``row expansion," which is
reflected in the left side of {\sevenrm (CAP)}.

\medbreak
{\bf Combinatorial Proof of Capelli's Identity}.
We will first figure out, step by step, what the combinatorial
objects that are being weight-enumerated by the left side of
{\sevenrm (CAP)}. Then we will decide who are the ``bad guys" and
will find an involution that preserves the absolute value of the
weight, but reverses the sign. The weight-enumerator of the
good guys will turn out to be counted by the right side of
{\sevenrm (CAP)}.

First we have to represent each $A_{ij}$ as a generating
polynomial over a particular set of combinatorial objects:
consider the 4-tuples $(i,j,k,l)$ where $i$, $j$,
$k=1,2,\ldots,n$ and $l=0,1$. For $i\not=j$ define 
${\cal A}_{ij}$ as the set of all 4-tuples $(a,b,c,d)$ such that
$a=i$, $c=j$, $d=0$ and $b=1,2,\ldots,n$. Next define
${\cal A}_{ii}$ as the set of all 4-tuples $(a,b,c,d)$ such
that $a=c=i$, and either $d=0$ and $b=1,2,\ldots,n$, or
$d=1$ and $b=i+1,\ldots,n$. Finally, let
$$
w(a,b,c,d)=\cases{x_{ba}\,p_{bc},&if $d=0$;\cr
h,&if $d=1$ (and $a=c$).\cr}
$$
We can then rewrite: $A_{ij}=\sum w(a,b,c,d)$, where
$(a,b,c,d)$ runs over all~${\cal A}_{ij}$. Hence
$$
\displaylines{(1.2)
\quad \smash{\sum_{\sigma \in {\cal S}_n}}\sgn (\sigma)\,
A_{\sigma 1,1}\ldots A_{\sigma n,n}\hfill\cr
\hfill{}=
\sum\sgn (a)\,
w(a_1,b_1,c_1,d_1)\ldots w(a_n,b_n,c_n,d_n),\quad\cr
}
$$
where the sum is over all sequences 
$(a_1,b_1,c_1,d_1,\ldots,a_n,b_n,c_n,d_n)$  satisfying
the properties: 
\item {1)} $a=(a_1,\ldots,a_n)$ is a permutation; 
\item {2)} $(c_1,\ldots,c_n)=(1,\ldots,n)$;
\item {3)} $d_i=0$ or 1 $(i=1,\ldots,n)$;
\item {4)} the $b_i$'s are arbitrary $(1\le b_i\le n)$ with the
sole condition that when $d_i=1$, then $a_i=i=c_i$ and 
$i+1\le b_i\le n$.

It then suffices to consider the set $\cal A$ of all
$4\times n$-matrices 
$$
G=\pmatrix{a_1&a_2&\ldots&a_n\cr
b_1&b_2&\ldots&b_n\cr
c_1&c_2&\ldots&c_n\cr
d_1&d_2&\ldots&d_n\cr}
$$
that satisfy the forementioned 1) to 4) properties 
and define the {\it weight} of $G$ as
$$
w(G)=\sgn (a) \prod_{i=1}^n
(x_{b_i,a_i}\,p_{b_i,i}(1-d_i)+hd_i).\leqno(1.3)
$$
Then the (1.2) sum may be expressed as:
$$
\sum_{\sigma \in {\cal S}_n}\sgn (\sigma) \,
A_{\sigma 1,1}\ldots A_{\sigma n,n}
=\sum_{G\in {\cal A}} w(G).
$$

If there is no pair $(i,j)$ such that $1\le i<j\le n$,
$d_i=d_j=0$ and $(b_i,i)=(b_j,a_j)$, say that $G$ is 
{\it not linkable}. Its weight can be expressed as a monomial
$\sgn(a )\,h^\alpha \prod x^\beta\prod p^\gamma$, where all the
$x$'s are written before all the $p$'s, by using the commutation
rule. If there exists such a pair $(b_i,i)=(b_j,a_j)$, the matrix
$G$ is said to be {\it linkable }. The product
$x_{b_i,a_i}p_{b_i,i}x_{b_j,a_j}p_{b_j,j}$  gives rise to the sum 
$x_{b_i,a_i}x_{b_i,i}p_{b_i,i}p_{b_j,j}
+x_{b_i,a_i}p_{b_j,j}h
=x_{b_i,a_i}x_{b_i,i}p_{b_i,i}p_{b_i,j}
+x_{b_i,a_i}p_{b_i,j}h$.
In the first monomial the commutation $x_{b_i,i}\,p_{b_i,i}$ has
been made; in the second monomial the latter product has
vanished and been replaced by~$h$. Such a pair  $(i,j)$ will be
called a {\it link}, of {\it source}~$i$ and {\it end}~$j$.

If a linkable matrix has $m$ links $(i_1<j_1)$, \dots, 
$(i_m<j_m)$, its weight will produce $2^m$ monomials when
the commutation rules are applied to it. Each of those
$2^m$ monomials corresponds to a subset
$K=\{{k_1},\ldots,{k_r}\}$ of the set $I=\{i_1,\ldots,i_m\}$ of
the link sources. We then have to consider the set of all the
pairs $(G,K)$, subject to the previous conditions and define the
weights of those pairs as single monomials in such a way that
the sum $\sum_K w(G,K)$ will be the weight of~$G$, once all the
commutations $px=xp+h$ have been made.

The weight $w(G,K)$ will be defined in the following way:
consider the single monomial introduced in (1.3); if $i$ belongs
to~$K$, drop $x_{b,i}$  and replace $p_{ b_i,i}$ by~$h$; if $i$
belongs to $I\setminus K$,  drop $x_{b,i}$ and replace 
$p_{b_i,i}$ by $x_{b_i,i}\,p_{b_i,i}$. Leave the other terms
alike.  In other words define the operators:
$$
\displaylines{
D_i=h {\partial\over \partial p_{b_i,i}}
{\partial\over \partial x_{b_i,i}}\qquad{\rm and}\qquad
\Delta_i=x_{b_i,i}p_{b_i,i}{\partial\over\partial p_{b_i,i}} 
{\partial\over\partial x_{b_i,i}}.\cr
\noalign{\hbox{Then  let}}
w(G,K)=\Bigl(\prod_{i\in K}D_i
\prod_{i\in I\setminus K}\Delta_i\Bigr) w(K)\cr} 
$$
For instance, the matrix
$$
G=\pmatrix{
4&5&1&8&7&6&9&2&3\cr
2&8&2&1&8&8&8&8&2\cr
1&2&3&4&5&6&7&8&9\cr
0&0&0&0&0&1&0&0&0\cr
}
$$
has three links $(1,3)$, $(2,8)$ and $(3,9)$. Its weight,
according to (1.3) reads:
$$-
x_{2,4}\,p_{2,1}\;
x_{8,5}\,p_{8,2}\;
x_{2,1}\,p_{2,3}\;
x_{1,8}\,p_{1,4}\;
x_{8,7}\,p_{8,5}\;
h\;
x_{8,9}\,p_{8,7}\;
x_{8,2}\,p_{8,8}\;
x_{2,3}\,p_{3,9}.
$$
Now consider the subset $K=\{2\}$ of its link source set
$I=\{1,2,3\}$. The weight of $w(G,K)$ is then:
$$-
x_{2,4}\,x_{2,1}\,p_{2,1}\;
x_{8,5}\,h\;
x_{2,3}\,p_{2,3}\;
x_{1,8}\,p_{1,4}\;
x_{8,7}\,p_{8,5}\;
h\;
x_{8,9}\,p_{8,7}\;
p_{8,8}\;
p_{3,9}.
$$
The simple drop-add rule just defined  guarantees that 
{\it no $p_{i,j}$ remains to the left of $x_{ij}$} in any of the
weight $w(G,K)$. After using all the commutations $px=xp+h$ we
then get 
$$
\sum_{\sigma \in {\cal S}_n} \sgn \sigma \,
A_{\sigma 1,1}\ldots A_{\sigma n,n}=\sum
w(G,K),\leqno(1.4) 
$$
where $G$ runs over all $\cal A$ and $K$ over all the
subsets of the link source set of~$G$.

It is obvious who the {\it good guys} are: those pairs $(G,K)$
such that  $G$ has no 1's on the last row and such that
$K$ is empty. The good guys correspond exactly to the
members of $\det X^t\,P$, in the classical case, where all
the $x_{i,j}$ commute with all the $p_{i,j}$, and obviously their
sum is  $\det X.\det P$. [A~combinatorial proof of which can be
found in [Z].] It remains to kill the bad guys, i.e., show that the
sum of their weights is zero.

If $(G,K)$ is a bad guy, $h$ occurs in $w(G,K)$ and either there
are 1's on the last row of~$G$, or $K$ is non empty. Let
$i=i(G,K)$ be the greatest integer $(1\le i\le n-1)$ such that
either $i$ a link source belonging to~$K$, or the $i$-th column 
has an entry equal to~1 on the last row.

In the first case, let $(i,j)$ be the link of source~$i$; then
replace the $i$-th  and the $j$-th columns as shown in the
next display, the other columns remaining intact:
$$
G=\pmatrix{&a_i&\ldots&i&\cr
&b_i&\ldots&b_i&\cr
&i&\ldots&j&\cr
&0&\ldots&0&\cr}\mapsto
\pmatrix{&i&\ldots&a_i&\cr
&j&\ldots&b_i&\cr
&i&\ldots&j&\cr
&1&\ldots&0&\cr}=G' .
$$
The link $(i,j)$ (with $i\in K$) has been suppressed. Let $(k,l)$ be
another link of~$G$ such that $k\in K$. Then $k<i$ by definition
of~$i$. On the other hand, $l\not=j$. If $l\not=i$, then $(k,l)$
remains a link of~$G'$. If $l=i$, then $(b_k,k)=(b_i,a_i)$ and the
link $(k,i)$ in~$G$ has been replaced by the link
$(k,j)$ in~$G'$, so that $k$ is still a link source in~$G'$.
Accordingly, $K\setminus\{i\}$ is a subset of the link set
of~$G'$ and it makes sense to define $K'=K\setminus\{i\}$. Also
notice that 
$$
i(G',K')=i(G,K).\leqno(1.5)
$$

As $G$ and $G'$ differ only by their $i$-th and $j$-th  columns,
the weights of $G$ and $G'$ will have {\it opposite sign};
furthermore, they will differ only by their $i$-th and $j$-th
factors, as indicated in the next display:  
$$
\vbox{\halign{\hfil$#\;$ &\hfil$\;#\;$\hfil
&\hfil$\;#\;$\hfil
&\hfil$\;#\;$\hfil
&\hfil$\;#\;$\hfil\cr
\left |w(G)\right |=\;\ldots 
&x_{b_i,a_i}\;p_{b_i,i}
&\ldots
&x_{b_i,i}\;p_{b_i,j}
&\ldots\cr
\left |w(G')\right |=\;\ldots
&h
&\ldots
&x_{b_i,a_i}\;p_{b_i,j}
&\ldots\cr
}}
$$
[The dots mean that the two words have the same left
factor, the same middle factor and the same right factor.]
Hence, as $i$ is in $K$, but not in $K'$, the operator
$D_i$ (resp. $\Delta_i$) is to be applied to $w(G)$ (resp.
$G'$) in order to get $w(G,K)$ (resp. $w(G,K')$), so that:
$$
\displaylines{
\vbox{\halign{\hfil$#\;$
&\hfil$\;#\;$\hfil
&\hfil$\;#\;$\hfil
&\hfil$\;#\;$\hfil
&\hfil$\;#\;$\hfil\cr
\left |w(G,K)\right |=\;\ldots 
&x_{b_i,a_i}\;h
&\ldots
&p_{b_i,j}
&\ldots\cr
\left |w(G',K')\right |=\;\ldots
&h
&\ldots
&x_{b_i,a_i}\;p_{b_i,j}
&\ldots\cr
}}\cr
\noalign{\hbox{showing that}}
\rlap{(1.6)}\hfill w(G,K)=-w(G',K').
\hfill\cr}
$$

In the second case the entries in the $i$-th column
$(i,j,i,1)$ satisfy the inequalities $i+1\le j\le n$, while
the $j$-th column (on the right of the $i$-th column) is
of the form $(a_j,b_j,j,0)$. Then define 
$$
G=\pmatrix{&i&\ldots&a_j&\cr
&j&\ldots&b_j&\cr
&i&\ldots&j&\cr
&1&\ldots&0&\cr}\mapsto
\pmatrix{&a_j&\ldots&i&\cr
&b_j&\ldots&b_j&\cr
&i&\ldots&j&\cr
&0&\ldots&0&\cr}=G',
$$
where only the $i$-th and $j$-th columns have been modified.
Clearly a new link $(i,j)$ has been created in $G'$. Let $(k,l)$ be
a link of~$G$ with $k\in K$. Then $k<i$. If $l=j$, we have
$(b_k,k)=(b_j,a_j)$, so that $(k,i)$ is a link of~$G'$. If
$l\not=j$, then $(k,l)$ remains a link of~$G'$. Thus $K\cup \{i\}$
is a set of link sources of~$G'$. It then makes sense to define
$K'=K\cup \{i\}$. Also notice that relation~(1.5) still holds.

As before, $w(G)$ and $w(G)$ have opposite signs. Furthermore 
$$
\displaylines{
\vbox{\halign{\hfil$#\;$
&\hfil$\;#\;$\hfil &\hfil$\;#\;$\hfil
&\hfil$\;#\;$\hfil
&\hfil$\;#\;$\hfil\cr
\left |w(G)\right |=\;\ldots 
&h
&\ldots
&x_{b_j,a_j}\;p_{b_j,j}
&\ldots\cr
\left |w(G')\right |=\;\ldots
&x_{b_j,a_j}\;p_{b_j,i}
&\ldots
&x_{b_j,i}\;p_{b_j,j}
&\ldots\cr
}}\cr
\noalign{\hbox{so that}}
\vbox{\halign{\hfil$#\;$
&\hfil$\;#\;$\hfil &\hfil$\;#\;$\hfil
&\hfil$\;#\;$\hfil
&\hfil$\;#\;$\hfil\cr
\left |w(G,K)\right |=\;\ldots 
&h
&\ldots
&x_{b_j,a_j}\;p_{b_j,j}
&\ldots\cr
\left |w(G',K')\right |=\;\ldots
&x_{b_j,a_j}\;h
&\ldots
&p_{b_j,j}
&\ldots\cr
}},\cr}
$$
showing that (1.6) also holds.

Taking into account (1.5) it is readily seen that
$\omega:(G,K)\mapsto (G',K')$ maps the first case into the second
one, and conversely. Applying $\omega$ twice gives the original
element, so it is an involution. Finally, property (1.6) makes it
possible to associate the bad guys into mutually canceling
pairs, and hence their total weight is zero.\qed

\medbreak
{\bf 2. A Combinatorial Proof of Turnbull's Identity.}

\proclaim Turnbull's Identity. Let $X=(x_{ij})$, $P=(p_{ij})$
$(1\le i,j\le n)$ be as before, but now they are symmetric
matrices: $x_{i,j}=x_{j,i}$ and $p_{i,j}=p_{j,i}$, their entries
satisfying the same commutation rules. Also let 
$\tilde P=(\tilde p_{i,j}):=(p_{i,j}(1+\delta _{i,j}))$. For each
positive integer~$n$ and for $1\le i,j\le n$, let
$$
\displaylines{\rlap{\rm (2.1)}\hfill
A_{ij}:=\sum_{k=1}^n
x_{ki}\tilde p_{kj} + h(n-i) \delta _{ij}.\hfill\cr
\noalign{\hbox{Then}}
\rlap{\sevenrm (TUR)}\hfill
\sum_{\sigma \in {\cal S}_n}
\sgn(\sigma )A_{\sigma 1,1}\ldots A_{\sigma n,n}=\det
X.\det \tilde P. \hfill\cr}
$$

The proof is very similar. However we have to introduce another
value for the $d_i$'s to account for the fact that the diagonal
terms of $\tilde P$ are $2p_{i,i}$. More precisely, for $i\not=j$
we let ${\cal T}_{i,j}$ be the set of all 4-tuples $(a,b,c,d)$
such that $a=i$, $c=j$, and either $d=0$ and $b=1,2,\ldots,n$, or
$d=2$ and $b=c=j$. In the same way, let ${\cal T}_{ii}$ be the
set of all 4-tuples $(a,b,c,d)$ such that $a=c=i$, and either
$d=0$ and $b=1,2,\ldots,n$, or $d=1$ and $b=i+1,\ldots,n$, or
$d=2$ and $b=c=i$. Finally, let 
$$ 
w(a,b,c,d)=\cases{x_{ba}\,p_{bc},&if $d=0$ or 2;\cr
h,&if $d=1$ (and $a=c$).\cr} 
$$ 
Next consider the set $\cal T$ of all $4\times n$-matrices 
$$
G=\pmatrix{a_1&a_2&\ldots&a_n\cr
b_1&b_2&\ldots&b_n\cr
c_1&c_2&\ldots&c_n\cr
d_1&d_2&\ldots&d_n\cr}
$$
satisfying the properties: 
\item{1)} $a=(a_1,\ldots,a_n)$ is a permutation; 
\item{2)} $(c_1,\ldots,c_n)=(1,\ldots,n)$;
\item{3)} $d_i=0$, 1 or 2 $(i=1,\ldots,n)$;
\item{ 4)} $b_i=\cases{1,\ldots,\ {\rm or}\ n,&when $d_i=0$;\cr 
i+1,\ldots,\ {\rm or}\ n,\ {\rm and}\ a_i=c_i=i,&when
$d_i=1$;\cr c_i=i,&when $d_i=2$.\cr}$

Then we have
$$
\sum_{\sigma \in {\cal S}_n}\sgn (\sigma) \,
A_{\sigma 1,1}\ldots A_{\sigma n,n}
=\sum_{G\in {\cal T}} w(G),
$$
where the weight $w(G)$ is defined as in (1.3) under the
restriction that the $d_i$'s are to be taken mod~2. 

Now to take the symmetry of $P$ and $X$ into account the
definition of a {\it link} has to be slightly modified. Say
that a pair $(i,j)$ is a {\it link} in~$G$, if $1\le i<j\le n$,
$d_i\equiv d_j\equiv 0\pmod 2$ and either
$(b_i,i)=(b_j,a_j)$, or $(b_i,i)=(a_j,b_j)$. In the Capelli case the
mapping $i\mapsto j$ (with $1\le i<j\le n$ and
$(b_i,i)=(b_j,a_j)$) set up a natural bijection of the source set
onto the end set. Furthermore, if the latter sets were of
cardinality~$m$, the weight of~$G$ gave rise to a polynomial
with $2^m$ terms. It is no longer the case in the Turnbull case.
For instance, if a matrix~$G$ is of the form
$$
G=\pmatrix{
\ldots&.&\ldots&.&\ldots&b_i&\ldots&i&\ldots\cr
\ldots&b_i&\ldots&i&\ldots&i&\ldots&b_i&\ldots\cr
\ldots&i&\ldots&b_i(=k)&\ldots&j&\ldots&l&\ldots\cr
\ldots&d_i&\ldots&d_k&\ldots&d_j&\ldots&d_l&\ldots\cr}
$$
with $d_i\equiv d_k\equiv d_j\equiv d_l\equiv 0\pmod 2$,
the weight of~$G$ will involve the factor
$$
p_{b_i,i}p_{i,b_i}x_{i,b_i}x_{b_i,i}=ppxx
$$
(by dropping the subscripts). The expansion of  the latter
monomial will yield
$$
ppxx=xxpp+4hxp+2h^2.
$$
With the term ``$xxpp$" all the commutations have been
made; say that {\it no link} remains. One link remains
unused to obtain each one of the next four terms ``$hxp$,"
i.e., $(i,j)$, $(i,l)$, $(k,j)$, $(k,l)$. 
Finally, the two pairs of links $\{(i,j),(k,l)\}$ and
$\{(i,l),(k,j)\}$ remain unused to produce the last term
``$2h^2$."

Accordingly, each of the term in the expansion of the
weight $w(G)$ (once all the commutations $px=xp+h$ have
been made) corresponds to a subset
$K=\{(i_1,j_1),\ldots,(i_r,j_r)\}$ of the link set of~$G$
having the property that all the $i_k$'s (resp. all the
$j_k$'s) are {\it distinct}. Let $w(G,K)$ denote the term
corresponding to $K$ in the expansion. We will then have
$$
w(G)=\sum _K w(G,K).
$$

As before the product $\det X.\det\tilde P$ is the sum
$\sum w(G,K)$ with $K$ empty and no entry equal to~1
on the last row of~$G$. If $(G,K)$ does not verify
the last two conditions, let $i=i(G,K)$ be the greatest
integer $(1\le i\le n-1)$ such that one of the following
conditions holds:

\item{1)} $d_i=0$  and $i$ is the  source of a link $(i,j)$
belonging to~$K$ such that $(b_i,i)=(b_j,a_j)$;

\item{2)} the $i$-th column  has an entry equal to~1 on the last
row;

\item{3)} $d_i=0$ or 2 and $i$ is the  source of a link $(i,j)$
belonging to~$K$ such that $(b_i,i)=(a_j,b_j)$ and case~1 does
not hold.

For cases 1 and 2 the involution $\omega: (G,K)\mapsto(G',K')$ 
is defined as follows:

\noindent
Case 1:
$$
G=\pmatrix{&a_i&\ldots&i&\cr
&b_i&\ldots&b_i&\cr
&i&\ldots&j&\cr
&0&\ldots&d_j&\cr}\mapsto
\pmatrix{&i&\ldots&a_i&\cr
&j&\ldots&b_i&\cr
&i&\ldots&j&\cr
&1&\ldots&d_j&\cr}=G' .
$$
Case 2:
$$
G=\pmatrix{&i&\ldots&a_j&\cr
&j&\ldots&b_j&\cr
&i&\ldots&j&\cr
&1&\ldots&d_j&\cr}\mapsto
\pmatrix{&a_j&\ldots&i&\cr
&b_j&\ldots&b_j&\cr
&i&\ldots&j&\cr
&0&\ldots&d_j&\cr}=G' .
$$
Notice that $d_j=0$ or 2 and when $d_j=2$, the matrix $G'$ also
belongs to~$\cal T$.

In case 1 the link $(i,j)$ has been suppressed. Let $(k,l)$ be a
link in~$G$ with $(k,l)\in K$. Then $k<i$ and $l\not=j$ because
of our definition of~$K$. If $l\not=i$, then $(k,l)$ remains a
link in~$G'$. Define $K'=K\setminus\{(i,j)\}$.

If $l=i$, then $(b_k,k)=(b_i,a_i)$ or
$(b_k,k)=(a_i,b_i)$ and the link $(k,i)$ in~$G$ has been
replaced by the link $(k,j)$ in~$G'$. In this case define
$K'=K\setminus \{(i,j),(k,i)\}\cup \{(k,j)\}$.
In those two subcases (1.5) remains valid.

In case 2 the link $(i,j)$ is now a link in~$G'$. Let $(k,l)$ belong
to $K$. Then $k<i$.  If $l\not=j$, then $(k,l)$ remains a link in
$G'$.  If $l=j$, then $(b_k,k)=(b_j,a_j)$ or ${}=(a_j,b_j)$, so
that $(k,i)$ is a link in~$G'$. Define $K'=K\cup\{(i,j)\}$ in the
first subcase and $K'=K\cup\{(i,j),(k,i)\}\setminus \{(k,j)\}$.
Again (1.5) holds.

If case 3 holds, $G$ has the form:
$$
G=\pmatrix{&a_i&\ldots&b_i&\cr
&b_i&\ldots&i&\cr
&i&\ldots&j&\cr
&d_i&\ldots&d_j&\cr}
$$
and eight subcases are to consider depending on whether $a_i$,
$b_i$ are equal or not to~$i$, and $d_i$ is equal to~0 or~2. The
two cases $a_i=b_i=i$ can be dropped, for $a$ is a permutation.
The two cases $b_i\not=i$, $d_i=2$ can also be dropped because
of condition~4 for the matrices in~$\cal T$. The case
$a_i\not=i$, $b_i=i$, $d_i=0$ is covered by case~1. There remain
three subcases for which the mapping $G\mapsto G'$ is defined
as follows:

\noindent
Case  $3'$: $a_i\not=i$, $b_i\not=i$, $d_i=0$.
$$
G=\pmatrix{&a_i&\ldots&b_i&\cr
&b_i&\ldots&i&\cr
&i&\ldots&j&\cr
&0&\ldots&d_j&\cr}\mapsto
\pmatrix{&b_i&\ldots&a_i&\cr
&a_i&\ldots&i&\cr
&i&\ldots&j&\cr
&0&\ldots&d_j&\cr}=G'.
$$
Case $3''$: $a_i=i$, $b_i\not=i$, $d_i=0$.
$$
G=\pmatrix{&i&\ldots&b_i&\cr
&b_i&\ldots&i&\cr
&i&\ldots&j&\cr
&0&\ldots&d_j&\cr}\mapsto
\pmatrix{&b_i&\ldots&i&\cr
&i&\ldots&i&\cr
&i&\ldots&j&\cr
&2&\ldots&d_j&\cr}=G'.
$$
Case $3'''$: $a_i\not=i$, $b_i=i$, $d_i=2$.
$$
G=\pmatrix{&a_i&\ldots&i&\cr
&i&\ldots&i&\cr
&i&\ldots&j&\cr
&2&\ldots&d_j&\cr}\mapsto
\pmatrix{&i&\ldots&a_i&\cr
&a_i&\ldots&i&\cr
&i&\ldots&j&\cr
&0&\ldots&d_j&\cr}=G'.
$$

In those three subcases the pair $(i,j)$ has remained a link
in~$G'$. Let $(k,l)$ be a link in~$K$ different from $(i,j)$. Then
$k<i$. Also $l\not=j$. If $l\not=i$, then $(k,l)$ remains a link
in $G'$. If $l=i$, then $(b_k,k)=(b_i,a_i)$ or $(b_k,k)=(a_i,b_i)$
and the link $(k,i)$ has been preserved in~$G'$. We can then
define: $K'=K$. Also (1.5) holds.

As for the proof of Capelli's identity we get
$w(G',K')=-w(G,K)$ in cases 1, 2 and~3. Clearly, $\omega$ maps
the first case to the second and conversely. Finally, subcase
$3'$ goes to itself, and $\omega$ exchanges the two subcases
$3''$ and $3'''$.

It follows that the sum of the weights of all the bad guys is
zero, thus establishing~{\sevenrm (TUR)}.

\medbreak
{\bf 3. What about the Anti-symmetric Analog?} Howe and Umeda
[H-U], and independently, Kostant and Sahi [K-S] discovered and
proved an anti-symmetric analog of Capelli's identity. Although
we, at present, are unable to give a combinatorial proof similar
to the above proofs, we state this identity in the hope that one
of our readers will supply such a proof. Since the
anti-symmetric analog is only valid for even~$n$, it is clear
that the involution cannot be ``local" as in the above
involutions, but must be ``global," i.e., involves many, if not all,
matrices.

\proclaim The Howe-Umeda-Kostant-Sahi Identity. Let $n$ be an
even positive integer. Let $X=(x_{i,j})$ $(1\le i,j\le n)$
be an anti-symmetric matrix: $x_{j,i}=-x_{i,j}$, and
$P=(p_{i,j})$ be the corresponding anti-symmetric momenta
matrix. Let
$$
\displaylines{
\rlap{$(1'')$}\hfill
A_{ij}:=\sum_{k=1}^n
x_{k,i}p_{k,j} + h(n-i-1) \delta _{ij}.\hfill\cr
\noalign{\hbox{Then}}
\rlap{\sevenrm (HU-KS)}\hfill
\sum_{\sigma \in {\cal S}_n}
\sgn(\sigma )A_{\sigma 1,1}\ldots A_{\sigma n,n}=\det
X.\det P. \hfill\cr}
$$

Although we are unable to prove the above identity
combinatorially, we do know how to prove combinatorially another,
less interesting, anti-symmetric analog of Capelli's identity,
that is stated without proof at the end of Turnbull's paper~[T].

\proclaim Turnbull's Anti-Symmetric Analog. Let $X=(x_{i,j})$ and
$P=(p_{i,j})$ $(1\le i,j\le n)$ be an anti-symmetric matrices as
above. Let
$$ 
\displaylines{
\rlap{$(1'')$}\hfill
A_{ij}:=\sum_{k=1}^n
x_{k,i}p_{k,j} - h(n-i) \delta _{ij},\hfill\cr
\noalign{\hbox{for $1\le i,j\le n$. Then}}
\rlap{\sevenrm (TUR$'$)}\hfill
\sum_{\sigma \in {\cal S}_n}
\sgn(\sigma )A_{\sigma 1,1}\ldots A_{\sigma n,n}={\rm
Per}(X^t\,P),\hfill\cr}
$$
where ${\rm Per}(A)$ denotes the permanent of a matrix~$A$, and
the matrix product $X^t\,P$ that appears on the right side of
{\sevenrm TUR$'$} is taken with the assumption that the
$x_{i,j}$ and $p_{i,j}$ commute.

Since the proof of this last identity is very similar to the proof
of Turnbull's symmetric analog (with a slight twist), we leave it
as an instructive and pleasant exercise for the reader.

\medskip
{\it Acknowledgement}. We should like to thank Roger
Howe for introducing us to Capelli's identity, and for
helpful conversations.

%reference for a paper
\def\article#1|#2|#3|#4|#5|#6|#7|
    {{\leftskip=7mm\noindent
     \hangindent=2mm\hangafter=1
     \llap{[#1]\hskip.35em}{#2}:
     #3, {\sl #4}, vol.\nobreak\ {\bf #5}, {\oldstyle #6},
     p.\nobreak\ #7.\par}}
%reference for a book
\def\livre#1|#2|#3|#4|
    {{\leftskip=7mm\noindent
    \hangindent=2mm\hangafter=1
    \llap{[#1]\hskip.35em}{#2}:
    {\sl #3}. #4.\par}}
%miscellaneous
\def\divers#1|#2|#3|
    {{\leftskip=7mm\noindent
    \hangindent=2mm\hangafter=1
     \llap{[#1]\hskip.35em}{#2}:
     #3.\par}}

\vskip 24pt

\centerline{REFERENCES}
\nobreak
\vskip 10pt
\article C|A. Capelli|\"Uber die Zur\"uckf\"uhrung der
Cayley'schen Operation $\Omega$ auf gew\"ohnliche
Polar-Operationen|Math. Annalen|29|1887|331-338|
\article G|L. Garding|Extension of a formula by Cayley to
symmetric determinants|Proc. Endinburgh Math. Soc. {\rm
Ser.~2}|8|1947|73--75|
\article H|R. Howe|Remarks on classical invariant
theory|Trans. Amer. Math. Soc.|313|1989|539--570|
\article H-U|R. Howe and T. Umeda|The Capelli identity,
the double commutant theorem, and multiplicity-free
actions|Math. Ann.|290|1991|565--619|
\article K-S|B. Kostant and S. Sahi|The Capelli Identity,
tube domains, and the generalized Laplace transform|Adv.
Math.|87|1991|71--92|
\article S|G. Shimura|On diferential operators attached
to certain representations of classical groups|Invent.
math.|77|1984|463--488|
\article T|H.W. Turnbull|Symmetric determinants and the
Cayley and Capelli operators|Proc. Edinburgh Math. Soc.
{\rm Ser.~2}|8|1947|73--75|
\livre T1|H.W. Turnbull|The Theory of Determinants,
Matrices, and Invariants|Dover, {\oldstyle 1960}|
\livre W|H. Weyl|The Classical Groups, their Invariants
and Representations|Princeton University Press,
{\oldstyle 1946}|
\article Z|D. Zeilberger|A combinatorial approach to
matrix algebra|Discrete Math.|56|1985|61--72|

\bye